\documentclass[hidelinks]{article}
\newcommand{\documentdate}{3 II 2026}

\usepackage{a4wide,latexsym,graphicx,amsmath}
\usepackage{varioref}
\usepackage{xcolor}
\usepackage{ulem}
\usepackage{amsfonts}
\usepackage{multirow}
\usepackage{hyperref}
\usepackage{todonotes}
\title{An objective-function-free algorithm for general smooth constrained optimization}
\author{Stefania Bellavia\footnotemark[1],
        Serge Gratton\footnotemark[2],
        Benedetta Morini\footnotemark[3],
        Philippe L. Toint\footnotemark[4]}

\newcommand{\beqn}[1]{\begin{equation}\label{#1}}
\newcommand{\eeqn}{\end{equation}}
\newcommand{\req}[1]{(\ref{#1})}
\newcommand{\ms}{\;\;\;\;}
\newcommand{\tim}[1]{\;\; \mbox{#1} \;\;}

\newcounter{algo}[section]
\renewcommand{\thealgo}{\thesection.\arabic{algo}}
\newcommand{\llem}[2]{\vspace{\baselineskip} 
\noindent\framebox[\textwidth]{\parbox{0.95\textwidth}{
\begin{lemma} \label{#1} \rm #2 \end{lemma} } } \vspace{\baselineskip} }
\newcommand{\algo}[3]{\refstepcounter{algo}
\begin{center}\begin{figure}[htbp]
\framebox[1.05\textwidth]{
\parbox{\textwidth} {\vspace{\topsep}
{\bf Algorithm \thealgo : #2}\label{#1}\\
\vspace*{-\topsep} \mbox{ }\\
{#3} \vspace{\topsep} }}
\end{figure}\end{center}}
\newcommand{\bpr}{{\bf Proof.} \hspace{1.5mm}}
\newcommand{\epr}{\hfill $\Box$ \vspace*{1em}}
\newcommand{\lthm}[2]{\vspace{\baselineskip} 
\noindent\framebox[\textwidth]{\parbox{0.95\textwidth}{
\begin{theorem} \label{#1} \rm #2 \end{theorem} } } \vspace{\baselineskip} }

\newcommand{\ii}[1]{\{ 1, \ldots, #1 \}}

\newcommand{\iibe}[2]{\{ #1, \ldots, #2 \}}

\newcommand{\calF}{{\cal F}} 
\newcommand{\calI}{{\cal I}} 
\newcommand{\calO}{{\cal O}}

\renewcommand{\Re}{\hbox{I\hskip -2pt R}}

\newcommand{\sfrac}[2]{{\scriptstyle \frac{#1}{#2}}}
\newcommand{\half}{\sfrac{1}{2}}
\newcommand{\eqdef}{\stackrel{\rm def}{=}}

\newcommand{\kap}[1]{\kappa_{\mbox{\tiny #1}}}

\newcommand{\flow}{f_{\rm low}}
\newcommand{\al}[1]{{\footnotesize{\sf #1}}}
\newcommand{\tal}[1]{{\normalsize {\sf #1}}}

\newcommand{\sT}[1]{s_{T,#1}}
\newcommand{\sN}[1]{s_{N,#1}}
\newcommand{\sTk}{\sT{k}}
\newcommand{\sNk}{\sN{k}}
\newcommand{\oT}[1]{\omega_{T,#1}}
\newcommand{\oN}[1]{\omega_{N,#1}}
\newcommand{\oTk}{\oT{k}}
\newcommand{\oNk}{\oN{k}}
\newcommand{\chiT}[1]{\chi_{T,#1}}
\newcommand{\chiN}[1]{\chi_{N,#1}}
\newcommand{\chiTk}{\chiT{k}}
\newcommand{\chiNk}{\chiN{k}}
\newcommand{\dT}[1]{d_{T,#1}}
\newcommand{\dN}[1]{d_{N,#1}}
\newcommand{\dTk}{\dT{k}}
\newcommand{\dNk}{\dN{k}}
\newcommand{\aT}[1]{\alpha_{T,#1}}

\newcommand{\aTk}{\aT{k}}

\newtheorem{theorem}{Theorem}[section]
\newtheorem{lemma}[theorem]{Lemma}
\newtheorem{corollary}[theorem]{Corollary}

\newcommand{\appnumsection}[1]{\section*{#1}
  \renewcommand{\theequation}{A.\arabic{equation}}
  \renewcommand{\thetheorem}{A.\arabic{theorem}}
  \renewcommand{\thetable}{A.\arabic{table}}
  \renewcommand{\thefigure}{A.\arabic{figure}}
  \renewcommand{\thesection}{A} }

\newcommand{\algname}{\al{ADIC}}
\newcommand{\proof}[1]{
\begin{list}{}{
\setlength{\topsep}{0.0pt}
\setlength{\partopsep}{0.0pt}
\setlength{\leftmargin}{0.025\textwidth}
\setlength{\rightmargin}{0.5\leftmargin}
\setlength{\labelwidth}{0.5\leftmargin}
\setlength{\labelsep}{0.25\leftmargin}}
\item \bpr #1 \epr \noindent
\end{list}}
\definecolor{checkgreen}{rgb}{0,0.6,0}
\newcommand{\comment}[1]{}

\begin{document}

\renewcommand{\thefootnote}{\fnsymbol{footnote}}
\maketitle

\footnotetext[1]{Dipartimento  di Ingegneria Industriale, Universit\`a degli Studi di Firenze,
  Firenze, Italia. Member of the INdAM Research Group GNCS.
  Email: stefania.bellavia@unifi.it.  Work partially supported by Progetti di Ricerca INDAM-GNCS.}
\footnotetext[2]{Universit\'e de Toulouse, INP, IRIT, Toulouse,
  France. Work partially supported by 3IA Artificial and Natural
  Intelligence Toulouse Institute (ANITI), French ``Investing for the
  Future - PIA3'' program under the Grant agreement ANR-19-PI3A-0004.
  Email: serge.gratton@toulouse-inp.fr.}
  \footnotetext[3]{Dipartimento  di Ingegneria Industriale, Universit\`a degli Studi di Firenze,
  Firenze, Italia. Member of the INdAM Research Group GNCS.
  Email: benedetta.morini@unifi.it.   Work partially supported by Progetti di Ricerca INDAM-GNCS.}
\footnotetext[4]{Namur Center for Complex Systems (naXys),
  University of Namur, Namur, Belgium.Work partially supported by 3IA Artificial and Natural
  Intelligence Toulouse Institute (ANITI) and DIEF (Florence).
  Email: philippe.toint@unamur.be.}
  
\renewcommand{\thefootnote}{\arabic{footnote}}

\begin{abstract}
A new algorithm for smooth constrained optimization is proposed that never computes the value of the problem's objective function  and that handles both equality and inequality constraints. The algorithm uses an adaptive switching strategy between a normal step aiming at reducing constraint's infeasibility and a tangential step improving dual optimality, the latter being inspired by the AdaGrad-norm method. Its worst-case iteration complexity is analyzed, showing that the norm of the gradients generated converges to zero like $\calO(1/\sqrt{k+1})$ for problems with full-rank Jacobians. Numerical experiments show that the algorithm's performance is remarkably
insensitive to noise in the objective function's gradient.
\end{abstract}

\hspace*{1.9cm} problems, reliability in the presence of noise, complexity.

\section{Introduction}\label{sec1}
The design and analysis of deterministic algorithms for solving constrained continuous optimization problems have a long history and have produced well-assessed techniques such as penalty methods, SQP methods, interior-point methods or filter methods (see \cite{NoceWrig99,ConnGoulToin00,GoulOrbaToin05,Yuan15} for example). These techniques all require the computation of both the function and derivative evaluation of the objective and the constraints.
By contrast, this paper addresses the solution of the problem 
\beqn{problem}
\min_{x \in \calF} f(x)
\tim{ where }
\calF = \{x\in\Re^n,\;\; c(x)=0\;\;  x\geq 0 \},
\eeqn 
using a first-oder objective function-free (OFFO) method.
Here $f$ is a smooth (possibly nonconvex) function from an open set containing the 
feasible region $\calF\subseteq \Re^n$ into $\Re$, $c(x): \Re^n \rightarrow \Re^m$
with $m\le n$ and equalities and inequalities are meant componentwise. Problem (\ref{problem}) includes general constrained optimization since all problems in this class can be 
cast into this form by using slack variables.  We assume that, given $x$, we can 
compute both the gradient $g(x) =\nabla_x f(x)$ of $f$ and the value of the 
constraints $c(x)$ as well as their Jacobian $J(x)=\nabla_x c(x)\in\Re^{m\times n}$,
which we will assume (for the purpose of our analysis) is full rank for any $x\ge 0$.

First-order OFFO procedures do not employ the value of the objective function but rely on gradients. They are known to be suitable for the solution of problems in which the function is approximated or subject to noise, and have exhibited remarkable robustness in the presence of noisy gradients \cite{BerCurtRobiZhou2021, CurtRobiZhou24,GratKopaToin23, GratToin25,WuWardBott18}.

Our approach is inspired by the "trust-funnel" approach \cite{GoulToin10,SampToin15}, which, as \cite{Zopp95}, has roots in the much older Himmelblau's "flexible tolerance " method \cite{Himm72}. As in these references, the new method uses an adaptive switching strategy to select a normal or a tangential steps, the first aiming at reducing the violation of the constraints and the second at improving the objective-function value. The stepsize in the latter is reminiscent of the stepsize formula used in Adagrad-norm method \cite{DuchHazaSing11,ward2020adagrad} for unconstrained optimization.
Our method can therefore be seen as a AdaGrad-like method for solving  equality and inequality  constrained problems. Other OFFO  methods for constrained optimization using AdaGrad stepsizes  have been proposed by the authors in the papers \cite{BellaGratMoriToin25, GratToin25}. In \cite{BellaGratMoriToin25} bound constrained optimization problems are considered; stochastic estimators of the gradient are allowed  and  second-order information are used when available.
The paper \cite{GratToin25}  considers nonlinear equality constrained problems with full-rank Jacobians and 
proposes a first-order algorithm that adaptively selects steps in the plane tangent to the constraints 
or steps that reduce infeasibility.  The  evaluation complexity is analyzed, in both cases yielding a 
global convergence rate in $O(1/\sqrt{k+1})$, identical in order to that of steepest-descent 
and Newton's methods for unconstrained problems \cite{CartGoulToin22}.

Our present proposal builds on these contributions and extends \cite {GratToin25} to handle inequality constraints and thus to cover general smooth constrained optimization. 
To accommodate such constraints, we revisited the procedure from \cite{GratToin25} by introducing suitable primal and dual criticality measures and redefining both tangential and normal steps, while avoiding a technical assumption on the first iteration.
Three different techniques are provided for the computation of the tangential step.   
We analyze the  worst-case iteration complexity of our procedures and show that the norm of the gradients generated converges to zero like $\calO(1/\sqrt{k+1})$.
Numerical experiments show that, in line with what happens on simpler problems, the algorithm's performance is remarkably
insensitive to noise in the objective function's gradient.

The authors are aware of four other papers on OFFO procedures \cite{BerCurtRobiZhou2021, CurtRobiZhou24, FangNaMahoKola24,Wangpierzhoucurt26} for constrained problems.  The paper 
\cite{BerCurtRobiZhou2021}
presents objective function-free Sequential Quadratic Programming (SQP) algorithms to solve smooth optimization problems with stochastic objective and  deterministic nonlinear equality constraints. It employs a stepsize selection scheme based on Lipschitz constants (or adaptively estimated Lipschitz constants) in place of the linesearch.  This approach has been extended in \cite{CurtRobiZhou24} to handle deterministic 
inequality constraints. A convergence analysis in expectation is carried out, but the worst-case complexity has not been analyzed. 
The method introduced in  \cite{FangNaMahoKola24} is designed to solve nonlinear optimization problems with stochastic objectives and deterministic equality constraints. It again
employs normal  and tangential steps, the latter being computed using  a standard trust-region technique;  an explicit penalty parameter is used and dynamically updated throughout the process, without requiring the objective function's computation.  Global almost-sure convergence is proved. \cite{Wangpierzhoucurt26} proposes a variant of the SQP approach of \cite{BerCurtRobiZhou2021} for equality-constrained problems with full-rank Jacobian using first-order methods with momentum and analyzes its rate of convergence. 

Our paper is organized as follows. The \tal{ADIC} (ADagrad with Inequality Constraints) class of algorithms is introduced in Section \ref{sec2} with its algorithmic options. Section \ref{sec3} analyzes its worst-case complexity. Results obtained from the numerical validation of the algorithms are described in Section\ref{sec4}. Section \ref{sec5} finally summarizes our contributions and discusses perspectives for further research.
\vskip 5pt \noindent
{\bf Notations:}  In what follows, $\|\cdot\|$ denotes the Euclidean norm unless otherwise specified, and $\sigma_{\min}[A]$ denotes the smallest singular value of the matrix $A$.

\section{The \tal{ADIC}\ class of algorithms}\label{sec2}

In the new class of iterative methods that we are going to describe, a
new iterate is formed using either a tangential step (i.e. a step in
the plane tangent to the manifold of equality constraints) or a normal
step (mostly orthogonal to that manifold),
the choice between the two being based on a comparison of the primal
and dual criticality measures.
One of the interesting features of this algorithmic framework is that it allows the use of fairly general non-negative
bounded dual and primal criticality measures,
denoted $\omega_T(x)$ and $\omega_N(x)$ respectively.

In the algorithm's description \vpageref{the-algo}, the successive iterates are denoted by $x_k$
and we let $g_k = g(x_k)$, $c_k=c(x_k)$, $J_k=J(x_k)$,
$\omega_T(x_k)=\oTk$, $\omega_N(x) = \oNk$.
These criticality measures are computed in Step~1, together with a
tangential stepsize whose form is, as we will detail later, directly
inspired by the AdaGrad \cite{DuchHazaSing11,ward2020adagrad} algorithm. Whether the step taken  is
tangential or normal is decided by comparing their sizes, each of
these steps being designed to provide a first-order improvement (of
a carefully chosen Lyapunov function) comparable to the relevant
criticality measure while being of a size ensuring that first-order
effects dominate (as we will prove below).

\algo{the-algo}{\algname($x_0$)}
{
\begin{description}
\item[Step 0: Initialization:] The constants
$\beta,\eta>0$, $\theta_T,\theta_N \ge 1$, $0<\eta_{\min}\le\eta_{\max}$ and $\kappa_t,\kappa_n, \varsigma \in (0,\half]$
 are given.\\
Project $x_0$ onto the positive orthant.\\
Set $k=0$ and $\Gamma_0 = 0$. 
\item[Step 1: Evaluations:]
  Evaluate $c_k=c(x_k)$, $J_k=J(x_k)$,  $g_k= \nabla f(x_k)$.
  Then compute the dual measure $\oT{k}$, the primal measure $\oNk$
  and the stepsize
  \beqn{aTk-long}
  \aTk = \frac{\eta}{\sqrt{\Gamma_k+\oTk^2+\varsigma}}.
  \eeqn
\item[Step 2: Normal  step: ] 
   Except possibly if
  \beqn{switch}
  \oN{k} \le \beta \,\aTk \oTk,
  \eeqn
  compute $\sNk$ such that
   \begin{eqnarray}
             & x_k+\sNk &\ge 0, \label{sN-feas}\\
             & \|\sNk\| &\le \theta_N \,\oNk, \label{sN-bound}
   \end{eqnarray}
   and there exists a contant $\kappa_n\in (0,\half)$ independent of $k$ such that
   \beqn{N-descent}
   \half\|c(x_k+\sNk)\|^2\le \half\|c_k\|^2 - \kappa_n\,\oNk^2.
   \eeqn
   Then set $x_k^+=x_k +\sNk$ .
   If \req{switch} holds and $\sNk$ was not computed, set $x_k^+ = x_k$.
\item[Step 3: Tangential  step: ]
   If \req{switch} holds, compute a step $\sTk$ such that
  \begin{eqnarray}
    x_k+\sTk & \ge & 0 \label{sT-feas}\\
    J_k \sTk & = & 0, \label{sT-nullspace}\\
    g_k^T\sTk & \le & -\kappa_t\,\aTk \,\oTk^2,\label{T-descent}\\
    \|\sTk\| & \le & \theta_T\,\aTk \,\oTk,\label{sT-bound}
    \end{eqnarray}
    and set $x_{k+1}=x_k^+ +\sTk$ and $\Gamma_{k+1}=\Gamma_k +\oTk^2$.\\
    Otherwise (i.e. if \req{switch} fails), set $x_{k+1}=x_k^+$ and $\Gamma_{k+1}=\Gamma_k$.

\item[Step 4: Loop: ] Increment $k$ by one and go to Step~1.
\end{description} 
}

\noindent
In our subsequent analysis, we need to  distinguish between iterates
using tangential or normal steps. We denote by
$\{k_\tau\}\subseteq \{k\}$ the index subsequence of iterations
such that a tangential step $\sTk$ was computed (implying that \req{switch} holds), while
$\{k_\nu\}$ is the index subsequence
of iterations where a normal step $\sNk$ was computed. Note that $\{k_\tau\}$ 
and $\{k_\nu\}$ need not be disjoint, but that $\{k\} = \{k_\tau\}\cup\{k_\nu\}$.
By convention, we will define $\sTk=0$ for $k\not\in\{k_\tau\}$ and
$\sNk = 0$ for $k\not\in \{k_\nu\}$.

We will also  consider the Lyapounov function (whose value is
hopefully decreased as the iterations progress) given by
\beqn{lyap-def}
\psi(x,\lambda) \eqdef L(x,\lambda)+\rho\|c(x)\|,
\eeqn
where $\rho$ is a fixed constant (to be determined below) and
$L(x,\lambda)$ is the standard Lagrangian 
\beqn{Lag-def}
L(x,\lambda) = f(x) + \lambda^Tc(x),
\eeqn
for some multiplier $\lambda\in \Re^m$.    The function $\psi(x,\lambda)$ is sometimes called the "sharp augmented Lagrangian" (see \cite{BuraGasiIsmaKaya06,BuraIuseMelo10,RomeFernTorr25} for instance). Of particular interest in our argument
is the least-squares Lagrange multiplier
$\widehat\lambda(x)$ defined by
\beqn{eq:ls-mult}
\big(J(x)J(x)^T\big)\,\widehat\lambda(x)\ =\ -\,J(x)\,g(x)
\eeqn
when the Jacobian $J(x)$ has full rank.

It is important to note to this point that, because all norms are
equivalent in $\Re^n$, our theoretically convenient choice of
expressing \req{sT-bound} and \req{sN-bound} in Euclidean norm is by
no means crucial.  Should other norms be used, as we will see below, the
relevant equivalence constants may be absorbed in $\theta_T$ and
$\theta_N$. It is also useful to notice \req{aTk-long} implies that
\beqn{achle1}
  \aTk \leq \frac{\eta}{\sqrt{\varsigma}}
  \tim{ and }
  \aTk \oTk < \eta.
\eeqn

Clearly, much else remains to be specified in our algorithmic outline:
details of which criticality measures are considered together with
which norm and methods to compute the
tangential step $\sTk$ itself as well as the normal step $\sNk$ must
be clarified. In order to simplify exposition, we focus in our theory
on a single technique for computing the normal step $\sNk$,
and propose to define it by one (or more) step(s) of a
trust-region algorithm applied on the constrained violation
$\half\|c(x)\|^2$ using a linear model. 
Lemma~\ref{TR-lemma} below will show that such a step satisfies our
requirements of Step~2  with
\beqn{chiNk-def}
\oN{k}=\chiNk = |c_k^TJ_k\dNk|.
\eeqn
where $\dNk$ solves
the problem
\beqn{LP-chiN}
   \begin{array}{cc}
      \min_d & c_k^T J_k d  \\
             & x_k+d \ge 0\\
             & \|d\|_\infty \le 1.\\
   \end{array}
\eeqn
By contrast, we will exploit the freedom in our model to introduce a
few variants for the computation of the tangential step.

\subsection{\tal{ADIC-LP}: two variants based on linear optimization}\label{ss:LP}

We start by describing a variant based on the dual criticality measure given by
\beqn{chiTk-def}
\oTk = \chiTk =|g_k^T\dTk|,
\eeqn
where $\dTk$ is the solution of the linear optimization\footnote{Formerly known as "linear programming".} problem
\beqn{LP-chiT}
   \begin{array}{cc}
      \min_d & g_k^Td  \\
        & J_k d=0\\
        & x_k+d \ge 0\\
        & \|d\|_\infty \le 1.
   \end{array}
\eeqn
Also observe that
\beqn{chiT-upper-LP}
\chiTk \le \|g_k\|\,\|\dTk\| \le \sqrt{n}\|g_k\|.
\eeqn
The tangential step $\sTk$ can then be computed in two ways. The first
is to define $\sTk$ as the solution of the linear 
programming problem
\beqn{LP-sT}
   \begin{array}{cc}
      \min_s  & g_k^Ts \\
        & J_k s=0\\
        & x_k+s \ge 0 \\
        & \|s\|_\infty\le \aTk \oTk.
   \end{array}
\eeqn 
(Note that \req{LP-sT} only differs from \req{LP-chiT} in the
definition of its bounds, and that we have used the liberty in the
choice of norms to express the bound on the step in $\|\cdot\|_\infty$).
A second, simpler, possibility is to
choose a multiple of $\dTk$ and simply set
\beqn{sT-backtrack}
\sTk = \frac{\aTk \oTk}{\|\dT{k}\|_\infty}\,\dT{k}.
\eeqn
Defined in either of these ways, $\sTk$ clearly satisfies
\req{sT-feas}, \req{sT-nullspace} and
\req{sT-bound} (with $\theta_T\ge \sqrt{n}$), and it is not difficult to verify that it also
satisfies \req{T-descent} with $\kappa_t=1/\max[\eta,1]$.

\llem{lemma:decrease-TS}{Suppose that, at tangential iteration
  $k_\tau$,  $\sT{k_\tau}$ is defined by either \req{LP-sT} or
  \req{sT-backtrack}. Then we have that 
\beqn{gen-decr}
|g_{k_\tau}^T\sT{k_\tau}|\ge \aT{k_\tau}\chiT{k_\tau}^2 = \frac{\aT{k_\tau}}{\max[\eta,1]}\oT{k_\tau}^2.
\eeqn
}

\proof{
Suppose first that $\|\sT{k_\tau}\|_\infty \ge \|\dT{k_\tau}\|_\infty$. Then $\dT{k_\tau}$ is feasible for problem \req{LP-sT} and thus
\[
|g_{k_\tau}^T\sT{k_\tau}| \ge |g_{k_\tau}^T\dT{k_\tau}| =\oT{k_\tau} 
\ge\frac{\aT{k_\tau}}{\eta} \oT{k_\tau}^2,
\]
where we used \req{achle1} to deduce the last inequality.
Suppose now that $\|\sT{k_\tau}\|_\infty < \|\dT{k_\tau}\|_\infty$. Then we must have that
$\|\sT{k_\tau}\|_\infty = \aT{k_\tau}\oT{k_\tau}$.
The vector 
$y = (\|\sT{k_\tau}\|_\infty/\|\dT{k_\tau}\|_\infty)\dT{k_\tau}$ is therefore feasible for problem \req{LP-chiT} and thus
\beqn{lin-dec}
|g_{k_\tau}^Ty| 
= \frac{\|\sT{k_\tau}\|_\infty}{\|\dT{k_\tau}\|_\infty}|g_{k_\tau}^T\dT{k_\tau}| 
= \frac{\aT{k_\tau}\oT{k_\tau}}{\|\dT{k_\tau}\|_\infty}\,\oT{k_\tau}
\ge \aT{k_\tau}\oT{k_\tau}^2
\eeqn
If $\sTk$ is defined by \req{sT-backtrack}, then $\sTk=y$ and
 \req{gen-decr} follows.  Otherwise, we obtain from \req{LP-sT}
 that $|g_{k_\tau}^T\sT{k_\tau}|\ge |g_{k_\tau}^Ty|$ and \req{gen-decr}
 also follows from \req{lin-dec}. 
}

\noindent
When the measure \req{chiTk-def} is used, it is also useful to note
that, for $x$ such that $c(x)=0$,
\beqn{chi-Lag}
\chi_{Tk}
=\chi_T(x_k)
=\min_{\|d\|_\infty\le 1} \{\nabla_x L(x_k,\widehat\lambda(x_k))^Td \mid J(x_k) d=0
      \tim{and} x_k+d \ge 0\} 
\eeqn
where the Lagrangian $L(x,\lambda)$ is defined in \req{Lag-def} and
$\widehat\lambda(x)$ is given by \req{eq:ls-mult}.

\subsection{\tal{ADIC-P1}: a projection-based variant}\label{ss:proj}

We next consider a variant based on the dual criticality measure given
by
\beqn{pik-def}
\oTk = \pi_T(x_k) \tim{ with }
\pi_T(x) = \| \Pi_\calF(x)\Big(x-g(x)\Big) - x\|
\eqdef \|p_1(x)\|.
\eeqn
where $\Pi_\calF(x)$ is the orthogonal projection onto
$\calF(x) \eqdef \{ x+y \in\Re^n \mid J(x)y = 0 \tim{and} x+y \ge 0 \}$
(see \cite[Section~12.1.4]{ConnGoulToin00}, for instance). In this
setting, one still defines $\aTk$  by \req{aTk-long} and one simply chooses
\beqn{proj1-sT}
\sTk = \min[\aTk,1] \, p_1(x_k).
\eeqn
Again, we note that
\beqn{chiT-upper-proj}
\pi_T(x_k) \le \|g_k\|.
\eeqn

The minimum in \req{proj1-sT} ensures that, 
by construction, $x_k+\sTk \in \calF(x_k)$ and thus that \req{sT-feas} holds.
The definition \req{proj1-sT} also implies that
\req{sT-nullspace} holds, while \req{sT-bound} with $\theta_T=1$ directly results from
\req{pik-def}.
The nature of the orthogonal projection also ensures the following
result.

\llem{lemma:decrease-proj}{Suppose that, at a tangential iteration
  $k_\tau$,  $\sT{k_\tau}$ is defined by \req{proj1-sT}. Then
\beqn{gen-decr-proj}
|g_{k_\tau}^T\sT{k_\tau}|
\ge \frac{\sqrt{\varsigma}}{\max[\eta,1]}\aT{k_\tau}\pi_T(x_{k_\tau})^2 
= \frac{\sqrt{\varsigma}}{\max[\eta,1]}\aT{k_\tau}\oT{k_\tau}^2.
\eeqn
}

\proof{The optimal nature of the projection implies that
  \[
  \Big([x_{k_\tau}-g_{k_\tau}]-[x_{k_\tau}+p_1(x_{k_\tau})]\Big)^T
  \Big([x_{k_\tau}+p_1(x_{k_\tau})]-x_{k_\tau}\Big) \ge 
  0
  \]
  and thus
  \[
  g_{k_\tau}^Tp_1(x_{k_\tau})
  \le -\|p_1(x_{k_\tau})\|^2
  = - \pi_T(x_{k_\tau})^2.
  \]
  Suppose first that $\sTk= \aTk p_1(x_k)$. Then, 
  \beqn{L22-case1}
  |g_{k_\tau}^T\sTk| \ge \aTk \pi_T(x_{k_\tau})^2.
  \eeqn
  Alternatively, if $\sTk= p_1(x_k)$, this implies that
  $\aTk \ge 1$. Now, \req{aTk-long} gives that 
  $\sqrt{\varsigma}\aTk \leq \max[\eta,1]$ and hence
  \beqn{L22-case2}
  |g_{k_\tau}^T\sTk| \ge \pi_T(x_{k_\tau})^2\ge\frac{\sqrt{\varsigma}}{\max[\eta,1]}\aTk\pi_T(x_{k_\tau})^2.
  \eeqn
  Combining \req{L22-case1} and \req{L22-case2} yields \req{gen-decr-proj}.
} 

\noindent
Thus the step \req{proj1-sT} also satisfies \req{T-descent} 
with $\kappa_t = \sqrt{\varsigma}/\max[\eta,1]$.
\subsection{Comments}

Some observations are in order at this stage.
\begin{enumerate}
\item Three types of iterations may occur in the course of the execution of the 
  algorithm. 
  \begin{itemize}
  \item The first is when the constraint violation is large, in which case 
        condition \req{switch} typically fails. A normal step $\sNk$ is then computed
        but a tangential step is not, which is probably reasonable because the meaning 
        of a  move in the tangent plane far away from the constraint is debatable, as it could result in very large steps which take forever to recover from.
  \item The second is when the constraint violation is moderate and \req{switch} holds.
        Both normal and tangential step may then be computed.
  \item The third is when the constraint violation is small. Condition 
        \req{switch} holds so that a tangential step is computed, but a normal step is not.
  \end{itemize} 
  What actually happens in a run depends on the choice of the constant $\beta$ in \req{switch} and the user's decision to avoid or force a normal step when possible.
\item The tangential stepsize formula \req{aTk-long} is
  of course reminiscent of the stepsize formula used in AdaGrad for
  unconstrained problems.
  Note that the running sum of squares of dual
  measures ($\Gamma_k$) is only updated at tangential iterations.
\item We have chosen to use the $\|.\|_\infty$ norm in \req{LP-chiN}, \req{LP-chiT} and
  \req{LP-sT}  so that these problems are standard
  linear programs, but, as we noted above, this is not necessary.  In particular,
  variants using the (isotropic) Euclidean norm or preconditioned
  version of these norms may also be considered.  One reason to
  consider Euclidean or other ellipsoidal norms is that
  $n$ inequality constraints created by the box constraints in the
  linear programs are replaced by a single constraint.
\item In our statement of the \algname\ framework, we have assumed
  that subproblems (\req{LP-chiT}, \req{LP-sT} or the projection
  problem in \req{pik-def}) are solved exactly. This is not necessary
  and it is sufficient that approximate solution are accurate enough to
  produce a decrease in the Lagrangian at least a fraction of the
  optimal one (as suggested by the introduction of the constant
  $\kappa_t$).
\end{enumerate}

\section{Worst-case complexity analysis}\label{sec3}
    
This section is devoted to the theoretical study of the
\algname\ method(s). Its main result is that, under suitable
conditions, the average value of $\left (\oTk+\|c_k\|\right)$ tends to
zero like $1/\sqrt{k+1}$. We need the following assumptions to derive it.
\begin{description}
\item[AS.0:] $f$ and $c$ are continuously differentiable on on open set containing the positive orthant of $\Re^n$.
\item[AS.1:] For all $x \ge 0$, $f(x)  \geq \flow$.
\item[AS.2:] For all $x\geq 0$, $\|g(x)\| \leq \kappa_g$ where $\kappa_g\ge \eta\beta$.
\item[AS.3:] For all $x\geq 0$, $\|c(x)\| \leq \kappa_c$, where $\kappa_c > 1$.
\item[AS.4:] For all $x\geq 0$, $\|J(x)\| \leq \kappa_J$
\item[AS.5:] For all $x\geq 0$, $\sigma_{\min}[J(x)] \geq \sigma_0 \in (0,1]$,
\item[AS.6:] The gradient $g(x)$ is globally Lipschitz continuous on the positive orthant (with constant $L_g$).
\item[AS.7:] The Jacobian $J(x)$ is globally Lipschitz continuous on the positive orthant
  (with constant $L_J$).
\item[AS.8:] There exists a constant $\xi\in(0,1]$ such that, for
  all $k\geq 0$, $\oNk \ge \xi \|c_k\|$.
\end{description}

Assumptions AS.1--AS.4  hold if the iterates remain, as is
often the case, in a closed bounded  set. Using the fact that the
product of bounded and Lipschitz functions is Lipschitz, we deduce the
following properties, whose detailed proofs can be found in appendix.

\llem{Lipschitz-things}{Suppose that AS.0 and AS.2--AS.7 hold. Then we have that
\begin{enumerate}
\item $c(x)$ is Lipschitz
  continuous on the positive orthant (with constant $L_c$),
\item $\nabla_x(\half \|c(x)\|^2) =
  J(x)^Tc(x)$ is Lipschitz continuous on the positive orthant (with constant $L_{JTc}\ge 1)$,
\item $\widehat\lambda(x)$ is well-defined on the positive orthant,
\item $\widehat\lambda(x)$ is bounded (by the constant $\kappa_\lambda$) and Lipschitz continuous (with constant $L_\lambda$) on the positive orthant,
\item $\nabla_xL(x,\lambda)$ is Lipschitz continuous on the positive orthant (with constant $L_L$).
\end{enumerate}
}

\noindent
We also observe  that \req{LP-chiN}, \req{chiNk-def}, AS.3 and AS.4 ensure that
\beqn{chiN-upper}
\chiNk \le \sqrt{n}\,\|J_k^Tc_k\|
\le \sqrt{n}\,\kappa_J\,\kappa_c,
\eeqn
Finally, AS.8 assumes that there exists a ``sufficient-descent''
direction for the problem \req{LP-chiN}. Specifically,  the normal
step is designed to reduce $\chiNk$ but it does not guarantee that
$\{\|c_{k_\nu}\|\}$ also converges to zero.  In fact, 
without further assumption,
the minimization of $\half\|c(x)\|^2$ may end up at a local minimizer
$x_{\rm loc}$ of this function which is infeasible for the original problem because
$c(x_{\rm loc})\neq 0$. The existence of such local minimizers may be
caused by a singular Jacobian $J(x_{\rm loc})$ (in which
case $J(x_{\rm loc})^Tc(x_{\rm loc})= 0$ does not imply $c(x_{\rm  loc}) = 0$),
or by the presence of bounds since $-J(x_{\rm loc})^Tc(x_{\rm loc})$
may then belong to the normal cone of the bound constraints at
$x_{\rm loc}$.
Unfortunately, convergence to such an $x_{\rm loc}$ cannot be avoided
without either applying a global optimization method to minimize
$\half\|c(x)\|^2$, or restricting the class of problems under
consideration.  Here we follow   the second approach and first  note that AS.5 already ensures
that $J(x_{\rm loc})^Tc(x_{\rm loc})= 0$  implies $c(x_{\rm  loc}) = 0$.
Making AS.8 is motivated by
the observation that descent along any direction for problem \req{LP-chiN}
not hitting the non-negativity constraints
must be limited by the bound $\|d\|_\infty\le 1$. Thus at least one
component of $\dNk$, say components $i\in\calI\subseteq \ii{n}$, must be equal to one in
absolute value, which implies that
$\chiNk = |c_k^TJ_k\dNk| \ge \|[J_k^Tc_k]_\calI\|_1$.
AS.8 then guarantees that $\|[J_k^Tc_k]_\calI\|_1$ is not negligible with
respect to $\|J_k^Tc_k\|\ge\sigma_0\|c_k\|$.

To maintain generality, we finally assume that the considered criticality measures are bounded.
\begin{description}
\item[AS.9:] There exists a constant $\kappa_\omega>0$ such that, for
  all $x\in \Re^n$, $\omega_T(x) \le \kappa_\omega$ and $\omega_N(x) \le \kappa_\omega$.
\end{description}
For the special cases discussed in Sections~\ref{ss:LP} and \ref{ss:proj}, AS.9
automatically results from AS.2--AS.4, as can be seen
from \req{chiT-upper-LP}, \req{chiT-upper-proj} and \req{chiN-upper}.

The next result shows our requirements on the normal step in Step~3 of
the \algname\ are not excessive. This is achieved by exhibiting one
particular computational scheme (a trust-region method) which satisfies
the conditions \req{sN-feas}--\req{N-descent}.

\llem{TR-lemma}{The normal step $\sNk$ (in Step~2) can be computed using a
  trust-region algorithm applied to minimizing  $\half\|c(x)\|^2$ subject 
  to \req{sN-feas} and \req{sN-bound} with $\|.\| =
  \|\cdot\|_\infty$ and $\oNk=\chiNk$ defined by \req{chiNk-def}, starting with the radius $\theta_N\chiNk$.
}
\proof{
  For any $\Delta \le \min[1,\theta_N\chiNk]$, let $s_{TR}(\Delta)$ be the solution of the
  problem of minimizing $c_k^TJ_ks$ over the constraints
  $x+s_{TR}(\Delta) \geq 0$ and $\|s_{TR}(\Delta)\|_\infty \le \Delta$.
  Using the Lipschitz continuity of $J(x)^Tc(x)$ (with constant
  $L_{JTc}$), we obtain that
  \beqn{TRl1}
  \half\Big(\|c(x_k+s_{TR}(\Delta))\|^2 - \|c_k\|^2\Big)
  \le c_k^TJ_ks_{TR}(\Delta) + \frac{L_{JTc}}{2}\Delta^2
  \le c_k^TJ_ks_{TR}(\Delta) + \frac{\max[\kappa_\omega,L_{JTc}]}{2}\Delta^2,
  \eeqn
  where $\kappa_\omega$ is defined in AS.9.
  Now, since $\Delta \le \min[1,\theta_N\chiNk]$, 
  given the vector $\dNk$ solution to \eqref{LP-chiN},   the vector $\Delta \dNk$ is feasible for
  \req{sN-feas}-\req{sN-bound}, and thus, from \req{chiNk-def},
  \[
  c_k^TJ_ks_{TR}(\Delta) \le c_k^TJ_k (\Delta \dNk) = -\chiNk\Delta.
  \]
  Hence
  \beqn{TRl2}
   c_k^TJ_ks_{TR}(\Delta)+ \frac{\max[\kappa_\omega,L_{JTc}]}{2}\Delta^2
  \le - \chiNk \Delta + \frac{\max[\kappa_\omega,L_{JTc}]}{2}\Delta^2.
  \eeqn
  It is then easy to verify that, if
  $\Delta \le \chiNk/\max[\kappa_\omega,L_{JTc}]$
  then \req{TRl2} gives that
  \beqn{TRl3}
  c_k^TJ_ks_{TR}(\Delta)+ \frac{\max[\kappa_\omega,L_{JTc}]}{2}\Delta^2
  \le -\frac{1}{2} \chiNk \Delta.
  \eeqn
 Remembering \req{chiN-upper}, we may then define
  \[
  \sNk= s_{TR}(\Delta_*)
  \tim{ with }
  \Delta_*
  = \min\left[ 1,\frac{\chiNk}{\max[\kappa_\omega,L_{JTc}]}\right]
  = \frac{\chiNk}{\max[\kappa_\omega,L_{JTc}]}
   \le \theta_N \chiNk
  \]
  where the last inequality, which shows that \req{sN-bound} holds, is derived using the bounds $L_{Tc} \ge 1$ and $\theta_N>1$.
Substituting this value in \req{TRl3} and using \req{TRl1} then yields that
  \[
  \half\Big(\|c(x_k+\sNk)\|^2 - \|c_k\|^2\Big)
  \le -\frac{\chiNk^2}{2\max[\kappa_\omega,L_{JTc}]}
  \]
  which proves the desired conclusion (with
  $\kappa_n=1/(2\max[\kappa_\omega,L_{JTc}]) \in (0,\half)$), because the radius $\Delta$
  can then be reduced (if necessary)  starting from $\theta_N\chiNk$ until
    \req{N-descent} holds.
} 

The previous result implies that the normal step can be computed using a
  trust-region algorithm for minimizing  $\half\|c(x)\|^2$ subject 
  to \req{sN-feas} and \req{sN-bound},  $\|.\| =
  \|\cdot\|_\infty$ and imposing that the trust-region solution $s_{TR}(\Delta)$ satisfies
  $$
   \half \|c(x_k+s_{TR}(\Delta))\|^2 \le \half \|c_k\|^2- \half \chiNk (\Delta)  .
  $$

Our analysis now proceeds by studying the behaviour of the Lyapunov function \req{lyap-def}
for iterations using normal and tangential steps (indexed by $k_\nu$
and $k_\tau$, respectively), before combining the results and
deriving global rates of convergence of $\left (\oTk+\|c_k\|\right)$
to zero along the sequences $\{k_\tau\}$,
$\{k_\nu\}$ and, finally, $\{k\}$.
For brevity, we define the abbreviated notations
\beqn{psi-def}
\psi(x)\eqdef \psi\big(x,\widehat\lambda(x)\big)
\tim{ and } \widehat\lambda_k = \widehat\lambda(x_k).
\eeqn
We also observe that  \req{Lag-def} and \req{eq:ls-mult} ensure that,
 for $\lambda^\dagger=\widehat\lambda(x)$,
\beqn{eq:ids}
\nabla_xL(x,\lambda^\dagger)= g(x)+J(x)^T \widehat\lambda(x)= g_T(x),
\eeqn
where $g_T(x)$ is the orthogonal projection of $g(x)$ onto the nullspace of $J(x)$,
and consequently, using AS.2,
\beqn{gradpsi}
\|\nabla_xL(x,\lambda^\dagger)\| \le \|g(x)\| \le \kappa_g.
\eeqn

\subsection{Descent at normal steps}

We first consider the effect of normal steps on the value of
the Lyapunov function $\psi$. We start by a very simple observation.

\llem{csq-decrease}{Suppose that AS.5 and AS.8 hold and that a normal step is used at
  iteration $k_\nu$.  Then $c_{k_\nu}^+=c(x_{k_\nu}+s_{N,k_\nu})$ satisfies
  \beqn{c-descent}
  \|c_{k_\nu}^+\| - \|c_{k_\nu}\| \le -\kappa_n\xi\,\oN{k_\nu} .
  \eeqn
}

\proof{
  We have from \req{N-descent} that
  $\|c_{k_\nu}^+\|<\|c_{k_\nu}\|$.  Then,
  \[
  2\|c_{k_\nu}\|(\|c_{k_\nu}\|-\|c_{k_\nu}^+\|)
  \ge (\|c_{k_\nu}\|+\|c_{k_\nu}^+\|)(\|c_{k_\nu}\|-\|c_{k_\nu}^+\|)
  = \|c_{k_\nu}\|^2-\|c_{k_\nu}^+\|^2,
  \]
  and therefore, using \req{N-descent} and AS.8, that
  \[
  \|c_{k_\nu}^+\|-\|c_{k_\nu}\| 
  \le -\frac{\kappa_n\oN{k_\nu}^2}{\|c_{k_\nu}\|}
  \le -\kappa_n \xi\, \oN{k_\nu}
  \]
} 

\noindent
We then use this observation to deduce the following result.

\llem{lem:normal-descent}{
Suppose that AS.3--AS.9 hold and that a normal step is used at
iteration $k_\nu$. Define
\beqn{rho-def}
\rho =  \frac{1}{\kappa_n\xi},\left[ (\kappa_g +\,\kappa_c\,\,L_\lambda)\theta_N\,+\,\left(\frac{\,L_L}{2}
  +\,L_\lambda L_c \right)\theta_N^2\kappa_\omega\,+\,\eta\;\right]
\eeqn
Then $x_{k_\nu}^+=x_{k_\nu}+s_{N,k_\nu}$ satisfies
\beqn{effectN}
\psi(x_{k_\nu}^+) - \psi(x_{k_\nu}) \leq -\eta \, \oN{k_\nu}.
\eeqn
}

\proof{We have that
\beqn{DaDb}
\psi(x_{k_\nu}^+) - \psi(x_{k_\nu})
= \underbrace{\psi(x_{k_\nu}^+,\widehat\lambda_{k_\nu})-\psi(x_{k_\nu},\widehat\lambda_{k_\nu})}_{\Delta_x}
+\underbrace{\psi(x_{k_\nu}^+,\widehat\lambda_{k_\nu}^+)-\psi(x_{k_\nu}^+,\widehat\lambda_{k_\nu})}_{\Delta_\lambda}.
\eeqn
Now consider $\Delta_x$ and $\Delta_\lambda$ separately. Using  the Lipschitz continuity of
$\nabla_x\psi(x,\widehat\lambda)$ ($\rho$ is fixed in \req{rho-def}) and
\req{c-descent}, we obtain that 
\beqn{dNe1}
\begin{aligned}
\Delta_x
& = \psi(x_{k_\nu}^+,\widehat\lambda_{k_\nu})-\psi(x_{k_\nu},\widehat\lambda_{k_\nu})\\
& = L(x_{k_\nu}^+,\widehat{\lambda}_{k_\nu})-L(x_{k_\nu},\widehat{\lambda}_{k_\nu})
+\rho\big(\|c_{k_\nu}^+\| - \|c_{k_\nu}\|\big)\\
& \leq (\nabla_xL(x_{k_\nu},\widehat\lambda_{k_\nu})^Ts_{N,k_\nu} +
r_3 - \rho \kappa_n\xi\,\oN{k_\nu}
\end{aligned}
\eeqn
with
\[
|r_3|\le \frac{L_L}{2}\|\sN{k_\nu}\|^2.
\]
We now invoke the Cauchy-Schwartz inequality, \req{gradpsi} and
\req{sN-bound} to deduce that
\beqn{eq:N1}
\begin{aligned}
\Delta_x
& \leq
\|\nabla_xL(x_{k_\nu},\widehat\lambda_{k_\nu})\|\,\|s_{N,k_\nu}\| -
\rho \frac{\kappa_n\xi}{2}\,\oN{k_\nu} + \frac{L_L}{2}\|\sN{k_\nu}\|^2\\
& \leq\kappa_g\|\sN{k_\nu}\|-  \rho \kappa_n\xi\,\oN{k_\nu} + \frac{L_L}{2}\|\sN{k_\nu}\|^2\\
& \leq \kappa_g\theta_N\oN{k_\nu}- \rho \kappa_n\xi\,\oN{k_\nu} + \frac{L_L}{2}\theta_N^2\oN{k_\nu}^2.
\end{aligned}
\eeqn
Using now the definition of 
$\Delta_\lambda$ in \req{DaDb}, AS.6, the Lipschitz continuity of $\widehat\lambda$ and $c$
and AS.3 then yields that
\beqn{eq:N2}
\begin{aligned}
\Delta_\lambda
& = \psi(x_{k_\nu}^+,\widehat\lambda(x_{k_\nu}^+))-\psi(x_{k_\nu}^+,\widehat\lambda(x_{k_\nu}))\\
&\le (\|c_{k_\nu}\|+\|c_{k_\nu}^+-c_{k_\nu}\|)\,\|\widehat\lambda_{k_\nu}^+- \widehat\lambda_{k_\nu}\|\\
&\le L_\lambda\,\|s_{N,k_\nu}\|\,\|c_{k_\nu}\|+L_\lambda L_c\|s_{N,k_\nu}\|^2\\
&\le L_\lambda\,\theta_N \kappa_c \oN{k_\nu}+L_\lambda  L_c\theta_N^2\oN{k_\nu}^2 \\
\end{aligned}
\eeqn
and thus, summing \req{eq:N1} and \req{eq:N2}, that
\[
\begin{aligned}
\psi(x_{k_\nu}^+) &- \psi(x_{k_\nu})\\
&\leq - \rho\kappa_n\xi\,\oN{k_\nu}+\kappa_g\theta_N\oN{k_\nu}
+L_\lambda\,\kappa_c\,\theta_N\oN{k_\nu}
+\left(\frac{\theta_N^2L_L}{2}+\theta_N^2 L_\lambda L_c \right)\oN{k_\nu}^2\\
&\leq - \rho\kappa_n\xi\oN{k_\nu}+\left(\kappa_g\theta_N
 +L_\lambda\,\kappa_c\,\theta_N\right)\oN{k_\nu}
 +\left(\frac{\theta_N^2L_L}{2}+ \theta_N^2L_\lambda L_c\right)\kappa_\omega\,\oN{k_\nu},
\end{aligned}
\]
where we have used AS.9 to deduce
the second inequality.
The bound \req{effectN} then follows from \req{rho-def}. \\
} 

\noindent
Note that, should $\sNk$ belong to the range space of $J_k$, the first term in the last right-hand side of \req{dNe1} vanishes and $\kappa_g$ disappears from \req{eq:N1} and, consequently, from \req{rho-def}.

\subsection{Descent at tangential steps}

We now turn to considering the effect of tangential steps.

\llem{tangent-decrease}{
Suppose that AS.4--AS.8 hold. Then
\beqn{effectT}
\psi(x_{k_\tau+1}) - \psi(x_{k_\tau}^+)
\le -\kappa_t\aT{k_\tau}\oT{k_\tau}^2
+\kap{tan}\,\aT{k_\tau}^2\oT{k_\tau}^2.
\eeqn
where
\beqn{kaptan-def}
\kap{tan} = \left[\frac{\theta_T^2}{2}\Big(L_L+\rho L_c\Big)
 + \beta\theta_N\theta_T\Big( L_L+\kappa_J L_\lambda+\rho L_J\Big)\right]
 + \frac{\beta\theta_TL_\lambda}{\xi}  + \theta_T^2 L_c L_\lambda.
\eeqn
}

\proof{
As in \req{DaDb}, we now have that
\beqn{DaDbtan}
\psi(x_{k_\tau+1}) - \psi(x_{k_\tau}^+)
= \underbrace{\psi(x_{k_\tau+1},\widehat\lambda_{k_\tau}^+)-\psi(x_{k_\tau}^+,\widehat\lambda_{k_\tau}^+)}_{\Delta_x}
+\underbrace{\psi(x_{k_\tau+1},\widehat\lambda_{k_\tau+1})-\psi(x_{k_\tau+1},\widehat\lambda_{k_\tau}^+)}_{\Delta_\lambda}.
 \eeqn
The Lipschitz
continuity of  $\nabla_x\psi(x,\widehat\lambda)$, \req{lyap-def} and \req{psi-def} give that
\beqn{Dx1}
\Delta_x = \nabla_x L(x_{k_\tau}^+,\widehat\lambda_{k_\tau}^+)^T \sT{k_\tau}
  + r_0 + \rho(\|c_{k_\tau+1}\|-\|c_{k_\tau}^+\|)
\tim{ with }
|r_0| \leq \frac{L_L}{2}\|\sT{k_\tau}\|^2.
\eeqn

{\color{blue}} Equation  (\ref{sT-nullspace}) gives 
\[
\begin{aligned}
\|c(x_{k_\tau+1})\|
& =\|c(x_{k_\tau}^+)-J_{k_\tau}\sT{k_\tau}+(J_{k_\tau}^+-J_{k_\tau}) \sT{k_\tau} + r_1\|\\
& \le \|c_{k_\tau}^+\| + \|r_1\| + \|J_{k_\tau}^+-J_{k_\tau}\|\,\|\sT{k_\tau}\|\\
& \le \|c_{k_\tau}^+\| + \|r_1\| + L_J\|\sN{k_\tau}\|\,\|\sT{k_\tau}\|\\
\end{aligned}
\]
with 
$\|r_1\| \leq \frac{L_c}{2} \|\sT{k_\tau}\|^2$.
Now $\|\sN{k_\tau}\|$ is either zero (if $k_\tau \not\in \{k_\nu\}$) 
or, using  \req{switch} for $k_\tau$, 
\[
\oN{k_\tau}\le  \beta \aT{k_\tau} \oT{k_\tau}
\]
and thus
\beqn{sNOsT}
\|\sN{k_\tau}\|
\le \theta_N\omega_{N,k_\tau}
\le \beta \theta_N \aT{k_\tau} \oT{k_\tau},
\eeqn
so that, whether $k_\tau \in \{k_\nu\}$ or not,  using \eqref{sT-bound},
\beqn{Encplus}
\|c(x_{k_\tau+1})\|-\|c(x_{k_\tau}^+)\|
\le \left(\frac{\theta_T^2L_c}{2}+\beta \theta_N\theta_TL_J\right)\aT{k_\tau}^2 \oT{k_\tau}^2.
\eeqn
Now differentiating $L$ with respect to its first argument and using
the Lipschitz continuity of $\nabla_xL$ with respect to this first
argument, AS.4, \eqref{sT-nullspace} and the Lipschitz continuity of $\widehat\lambda$ gives that
\[
\begin{aligned}
\nabla_x L(x_{k_\tau}^+,\widehat\lambda_{k_\tau}^+)^T \sT{k_\tau}
&= \Big(\nabla_x L(x_{k_\tau}^+,\widehat\lambda_{k_\tau}^+)^T \sT{k_\tau}-\nabla_x L(x_{k_\tau},\widehat\lambda_{k_\tau}^+)^T \sT{k_\tau}\Big)\\
& \hspace*{1cm} +\Big(\nabla_xL(x_{k_\tau},\widehat\lambda_{k_\tau}^+)^T\sT{k_\tau} - \nabla_x L(x_{k_\tau},\widehat\lambda_{k_\tau})^T\sT{k_\tau}\Big)\\
& \hspace*{1cm} +g_{k_\tau}^T\sT{k_\tau} + \widehat\lambda_{k_\tau}^TJ_{k_\tau}\sT{k_\tau}\\
&= \Big(\nabla_x L(x_{k_\tau}^+,\widehat\lambda_{k_\tau}^+)^T \sT{k_\tau}-\nabla_x L(x_{k_\tau},\widehat\lambda_{k_\tau}^+)^T \sT{k_\tau}\Big)\\
& \hspace*{1cm} +\Big( (\widehat\lambda_{k_\tau}^+)^TJ_{k_\tau}  - \widehat\lambda_{k_\tau}^TJ_{k_\tau}\Big)^T\sT{k_\tau}+g_{k_\tau}^T\sT{k_\tau} + \widehat\lambda_{k_\tau}^TJ_{k_\tau}\sT{k_\tau}\\
&\le \Big( L_L+\kappa_J L_\lambda\Big)\|\sN{k_\tau}\|\,\|\sT{k_\tau}\| +g_{k_\tau}^T\sT{k_\tau} + \widehat\lambda_{k_\tau}^TJ_{k_\tau}\sT{k_\tau}\\
&\le g_{k_\tau}^T\sT{k_\tau} + \beta\Big( L_L+\kappa_J L_\lambda\Big)\theta_N\theta_T\aT{k_\tau}^2 \oT{k_\tau}^2\\
\end{aligned}
\]
where the last inequality results from \req{sT-nullspace} and \req{sNOsT}.
Hence we obtain from \req{Dx1}, \req{T-descent}, \req{sT-bound}  and \req{Encplus}  that
\beqn{desc-det}
\begin{aligned}
\Delta_x  
&\le g_{k_\tau}^T\sT{k_\tau} + r_0
+ \rho\left(\frac{\theta_T^2L_c}{2}+\beta\theta_N\theta_TL_J\right)\aT{k_\tau}^2 \oT{k_\tau}^2 + \beta \theta_N\theta_T\Big( L_L+\kappa_J L_\lambda\Big)\,\aT{k_\tau}^2 \oT{k_\tau}^2 \\
&\le g_{k_\tau}^T\sT{k_\tau}
+\left(\frac{\theta_T^2}{2}\left(L_L+\rho L_c\right)+\beta\theta_N\theta_T\Big( L_L+\kappa_J L_\lambda+\rho L_J\Big)\right)\,\aT{k_\tau}^2 \oT{k_\tau}^2 \\
&\le -\kappa_t\aT{k_\tau}\oT{k_\tau}^2
+\left(\frac{\theta_T^2}{2}\left(L_L+\rho L_c\right)+\beta\theta_N\theta_T\Big( L_L+\kappa_J L_\lambda+\rho L_J\Big)\right)\,\aT{k_\tau}^2 \oT{k_\tau}^2.
\end{aligned}
\eeqn
Now, we may use the Lipschitz continuity of $\widehat\lambda$ and
$c$, inequality \req{N-descent}, the Cauchy-Schwartz inequality,  and AS.8 to deduce that
\beqn{Dlprod}
\begin{aligned}
\Delta_\lambda
& =c_{k_\tau+1}^T\big(\widehat\lambda_{k_\tau+1}-\widehat\lambda_{k_\tau}^+\big)\\
&= (c_{k_\tau+1}-c_{k_\tau}^+)^T\big(\widehat\lambda_{k_\tau+1}-\widehat\lambda_{k_\tau}^+\big)
 + (c_{k_\tau}^+)^T\big(\widehat\lambda_{k_\tau+1}-\widehat\lambda_{k_\tau}^+\big)\\
&\le \|c_{k_\tau}^+\|\,\|\widehat\lambda_{k_\tau+1}-\widehat\lambda_{k_\tau}^+\|
 + \|c_{k_\tau+1}-c_{k_\tau}^+\|\,\|\widehat\lambda_{k_\tau+1}-\widehat\lambda_{k_\tau}^+\|\\
&\le \|c_{k_\tau}\|\,\|\widehat\lambda_{k_\tau+1}-\widehat\lambda_{k_\tau}^+\|
 + \|c_{k_\tau+1}-c_{k_\tau}^+\|\,\|\widehat\lambda_{k_\tau+1}-\widehat\lambda_{k_\tau}^+\|\\
&\le \frac{L_\lambda}{\xi} \oN{k_\tau} \|\sT{k_\tau}\| + L_c L_\lambda\|\sT{k_\tau}\|^2.
\end{aligned}
\eeqn
Again using  \req{switch} for $k\in\{k_\tau\}$, \req{sT-bound} and \req{sNOsT},
we obtain that
\[
\Delta_\lambda
\le \frac{\beta\theta_T L_\lambda}{\xi}\aT{k_\tau}^2 \oT{k_\tau}^2 +
    \theta_T^2L_c L_\lambda\aT{k_\tau}^2 \oT{k_\tau}^2.
\]
Thus, summing $\Delta_x$ and $\Delta_\lambda$, we deduce that
\[
\begin{aligned}
\psi(x_{k_\tau+1}) - \psi(x_{k_\tau}^+)
& \leq -\kappa_t\aT{k_\tau} \oT{k_\tau}^2
+\left(\frac{\theta_T^2}{2}\Big(L_L+\rho L_c\Big)+\beta\theta_N\theta_T\Big( L_L+\kappa_J L_\lambda+\rho L_J\Big)\right)\aT{k_\tau}^2\oT{k_\tau}^2\\
&\hspace*{1cm}+\left(\frac{\beta\theta_T L_\lambda}{\xi} +\theta_T^2 L_c L_\lambda\right)\aT{k_\tau}^2\oT{k_\tau}^2\\
\end{aligned}
\]
and \req{effectT} follows.
} 

\noindent
Observe that the second term in the bracket of \req{kaptan-def} only appears when $k_\tau\in\{k_\nu\}$.
The bound \req{effectT} quantifies the effect of  tangential steps
on the Lyapunov function, and its right-hand side involves a
first-order (descent) term and a second-order perturbation term.
We now derive crucial bounds on these terms, using the fact that
$\Gamma_k$ is not updated at normal iterations.

\llem{lem:adagrad}{
Suppose that AS.2 and AS.5 hold. If we denote
\[
\Gamma_{k_{\tau_0}} = 0, 
\quad
\Gamma_{k_{\tau+1}} =\Gamma_{k_{\tau}}+\oT{k_\tau}^2,
\quad
\aT{k_\tau} =\frac{\eta}{\sqrt{\varsigma+\Gamma_{k_{\tau+1}}}},\quad
\]
then, for all $\tau_0\le \tau_1$,
\beqn{eq:AG1}
\sum_{\tau=\tau_0}^{\tau_1}\aT{k_\tau}\,\oT{k_\tau}^2
> \eta\sqrt{\varsigma}\,\sqrt{1+\frac{\Gamma_{k_{\tau_1+1}}}{\varsigma}} - \eta\sqrt{\varsigma}
\eeqn
\beqn{eq:AG2}
\sum_{\tau=\tau_0}^{\tau_1}\aT{k_\tau}^2\,\oT{k_\tau}^2
\le
\eta^2\,\log\left(1+\frac{\Gamma_{k_{\tau_1+1}}}{\varsigma}\right).
\eeqn
}

\proof{
Let $w_{k_\tau+1} =\sqrt{\Gamma_{k_{\tau+1}}+\varsigma}$.
The definition of $\aT{k_\tau}$ in \req{aTk-long} implies that
\[
\begin{aligned}
\sum_{\tau=\tau_0}^{\tau_1} \aT{k_\tau} \,\oT{k_\tau}^2
& =    \eta \sum_{\tau=\tau_0}^{\tau_1}
\frac{\oT{k_\tau}^2}{\sqrt{\varsigma+\Gamma_{k_{\tau+1}}}}\\
&> \eta \sum_{\tau=\tau_0}^{\tau_1} \frac{\oT{k_\tau}^2}{w_{k_{\tau+1}}+w_{k_{\tau}}}\\
&= \eta \sum_{\tau=\tau_0}^{\tau_1} \frac{w_{k_{\tau+1}}^2-w_{k_{\tau}}^2}{w_{k_{\tau+1}}+w_{k_{\tau}}}\\
&= \eta \sum_{\tau=\tau_0}^{\tau_1} (w_{k_{\tau+1}}-w_{k_{\tau}})\\
&= \eta\big(w_{k_{\tau_1+1}}-w_{k_{\tau_0}}\big).
\end{aligned}
\]
Now observe that, using $\Gamma_{k_{\tau_0}} = 0$,
\[
w_{k_{\tau_1+1}}-w_{k_{\tau_0}}
= \sqrt{\varsigma+\Gamma_{k_{\tau_1+1}}}-\sqrt{\varsigma+\Gamma_{k_{\tau_0}}}\\
=\sqrt{\varsigma+\Gamma_{k_{\tau_1+1}}}-\sqrt{\varsigma},
\]
which then gives \req{eq:AG1}.  Using the concavity and the increasing
nature of the logarithm, we also
have from \req{aTk-long} that
\[
\aT{k_\tau}^2\oT{k_\tau}^2\!
= \eta^2  \frac{\oT{k_\tau}^2}{\varsigma+\Gamma_{k_\tau+1}}\!
= \eta^2  \frac{\Gamma_{k_{\tau+1}}-\Gamma_{k_{\tau}}}{\varsigma+\Gamma_{k_{\tau+1}}}
\leq \eta^2 \left[
\log(\varsigma+\Gamma_{k_{\tau+1}})-\log(\varsigma+\Gamma_{k_{\tau}})\right].
\] 
Summing for $\tau\in\{\tau_0, \ldots, \tau_1\}$ then yields that
\[
\sum_{\tau=\tau_0}^{\tau_1}\aT{k_\tau}^2\,\oT{k_\tau}^2
\le \eta^2\left[\log(\varsigma+\Gamma_{k_{\tau_1+1}})-\log(\varsigma+\Gamma_{k_{\tau_0}})\right],
\]
and \req{eq:AG2} follows, again using $\Gamma_{k_{\tau_0}} = 0$.
} 

\subsection{Telescoping sum}

Having considered the impacts of tangential and normal steps
separately, we now combine them to derive a crucial inequality.

\llem{telescoping}{
Suppose that AS.0--AS.8 hold. Then, for any $\tau_1 > 0$ and any
$\nu_1 \ge 0$,
\beqn{full-telescopic}
\sqrt{1+\frac{\Gamma_{k_{\tau_1+1}}}{\varsigma}} + \sum_{\nu=\nu_0}^{\nu_1}\oN{k_\nu}
\leq \kap{gap} + \frac{\kap{tan}}{\kappa_t\sqrt{\varsigma}}\log\left(1+\frac{\Gamma_{k_{\tau_1+1}}}{\varsigma}\right),
\eeqn
where
\[
\kap{gap} = \frac{1}{\eta\kappa_t\sqrt{\varsigma}}\left(1+\psi(x_0) +
\kappa_c\kappa_\lambda +\rho\kappa_c - \flow\right).
\]
}

\proof{
Consider $k \ge 0$. Then, defining
$\min[k_{\nu_0},k_{\tau_0}] = 0$ and $\max[k_{\nu_1},k_{\tau_1}]= k$,
we have that $\Gamma_{k_{\tau_0}}=0$ and we may apply Lemma~\ref{lem:adagrad}.
Combining \req{effectT}, \req{eq:AG1} and \req{eq:AG2}, we obtain that
\beqn{eq:sumTtg}
\begin{aligned}
\sum_{\tau=\tau_0}^{\tau_1}\left(\psi(x_{k_\tau+1})-\psi(x_{k_\tau}^+)\right)
&\le \eta\kappa_t\sqrt{\varsigma}-\eta\kappa_t\sqrt{\varsigma}\,
\sqrt{1+\frac{\Gamma_{k_{\tau_1+1}}}{\varsigma}}
 +\eta^2\kap{tan}\log\left(1+\frac{\Gamma_{k_{\tau_1+1}}}{\varsigma}\right).
\end{aligned}
\eeqn
Also considering \req{effectN} and observing that
$x_k^+ = x_{k+1}$ when $k\in\{k_\nu\}\setminus\{k_\tau\}$
and $x_k^+ = x_k$ when $k\in\{k_\tau\}\setminus\{k_\nu\}$
therefore yields that
\beqn{eq:sumTtg2}
\begin{aligned}
\psi(x_{k+1})- \psi(x_0)
&=\sum_{\tau=\tau_0}^{\tau_1}\big(\psi(x_{k_\tau+1})-\psi(x_{k_\tau}^+)\big)
+ \sum_{\nu=\nu_0}^{\nu_1}\big(\psi(x_{k_\nu}^+)-\psi(x_{k_\nu})\big) \\
& \le \eta\kappa_t\sqrt{\varsigma} - \eta\kappa_t\sqrt{\varsigma}\,\sqrt{1+\frac{\Gamma_{k_{\tau_1+1}}}{\varsigma}}
      - \eta \sum_{\nu=\nu_0}^{\nu_1}\oN{k_\nu}
      + \eta^2\kap{tan}\log\left(1+\frac{\Gamma_{k_{\tau_1+1}}}{\varsigma}\right)\\
& \le \eta\kappa_t\sqrt{\varsigma}- \eta\kappa_t\sqrt{\varsigma}\,\sqrt{1+\frac{\Gamma_{k_{\tau_1+1}}}{\varsigma}}
      - \eta \kappa_t\sqrt{\varsigma}\,\sum_{\nu=\nu_0}^{\nu_1}\oN{k_\nu}
+ \eta\kap{tan}\log\left(1+\frac{\Gamma_{k_{\tau_1+1}}}{\varsigma}\right),
\end{aligned}
\eeqn
where we used the facts that $\kappa_t\sqrt{\varsigma}<1$ and $\eta\le 1$.
Using now \req{lyap-def}, \req{psi-def}, the Cauchy-Schwartz inequality, the boundedness of $\widehat\lambda(x)$, AS.1 and AS.3, we have that
\[
\psi(x_{k+1}) -\psi(x_0) -\eta\kappa_t\sqrt{\varsigma}\geq \flow- \kappa_c\kappa_\lambda-\rho\kappa_c-\psi(x_0)-\eta\kappa_t\sqrt{\varsigma}
\eqdef -\eta\kappa_t\sqrt{\varsigma}\,\kap{gap},
\]
so that \req{eq:sumTtg2} implies \req{full-telescopic}.
} 

\subsection{Tangential complexity}

Lemma~\ref{telescoping} implies upper bounds for both $\Gamma_{k\tau}$ and
$\oN{k_\nu}$. We now exploit the first of these to derive the rate of
convergence for tangential steps proper, after establishing a useful
technical result. 

\llem{tech-tan}{
  Suppose that $at \le b + c\,\log(t)$  for $t\geq 1$ and $a,c>0$.  Then, 
  \[
  t \le \frac{2b}{a}+\frac{2c}{a}\left[\log\left(\frac{2c}{a}\right)-1 \right].
  \]
}

\proof{ See \cite[Lemma~3.2]{BellaGratMoriToin25}.
} 

\llem{tangent-complexity}{
Suppose that AS.0--AS.9 hold.  Then, for any $\tau_1> 0$,
\beqn{Gamma-bound}
\sqrt{\varsigma+\Gamma_{k_{\tau_1+1}}}\le \kappa_T
\eqdef
2 \kap{gap}\sqrt{\varsigma} + \frac{4\kap{tan}}{\kappa_t}\left[\log\left(\frac{4\kap{tan}}{\kappa_t\sqrt{\varsigma}}\right)-1\right]
\eeqn
and
\beqn{Tcomp}
\xi\sum_{\tau=\tau_0}^{\tau_1}\Big(\oT{k_\tau}+\|c_{k_\tau}\|\Big)
\leq \sum_{\tau=\tau_0}^{\tau_1}\Big(\oT{k_\tau}+\oN{k_\tau}\Big)
\leq \kappa_T\sqrt{\tau_1+1}\left(1+\frac{\beta\eta}{\sqrt{\varsigma}}\right).
\eeqn
}

\proof{
The bound \req{full-telescopic} implies that
\[
\sqrt{1+\frac{\Gamma_{k_{\tau_1+1}}}{\varsigma}}
\leq \kap{gap} + \frac{\kap{tan}}{\kappa_t\sqrt{\varsigma}}\log\left(1+\frac{\Gamma_{k_{\tau_1+1}}}{\varsigma}\right)
= \kap{gap} + \frac{2\kap{tan}}{\kappa_t\sqrt{\varsigma}}\log\left(\sqrt{1+\frac{\Gamma_{k_{\tau_1+1}}}{\varsigma}}\right).
\]
Using Lemma~\ref{tech-tan} with
\[
t = \sqrt{1+\frac{\Gamma_{k_{\tau_1+1}}}{\varsigma}},
\ms\ms
a = 1,
\ms\ms
b = \kap{gap}
\ms\tim{and }\ms
c = \frac{2\kap{tan}}{\kappa_t\sqrt{\varsigma}},
\]
we then obtain that
\[
\sqrt{\varsigma+\Gamma_{k_{\tau_1+1}}}
= \sqrt{\varsigma}\,\sqrt{1+\frac{\Gamma_{k_{\tau_1+1}}}{\varsigma}}
\le \sqrt{\varsigma}\left\{2\kap{gap}+ \frac{4\kap{tan}}{\kappa_t\sqrt{\varsigma}}
\left[\log\left(\frac{4\kap{tan}}{\kappa_t\sqrt{\varsigma}}\right)-1\right]\right\}. 
\]
This is \req{Gamma-bound}.
We may now invoke the inequality
\[
\sum_{j=0}^k a_j \le \sqrt{k+1}\sqrt{\sum_{j=0}^k a_j^2}
\]
for nonnegative $\{a_j\}_{j=0}^k$ to deduce from the definition of
$\Gamma_{k_\tau}$ and \req{Gamma-bound} that
\beqn{tantan}
\sum_{\tau=\tau_0}^{\tau_1} \oT{k_\tau}
\le \sqrt{\tau_1+1}\,\sqrt{\sum_{\tau=\tau_0}^{\tau_1}\oT{k_\tau}^2}
= \sqrt{\tau_1+1}\,\sqrt{\Gamma_{k_{\tau_1}+1}}
< \sqrt{\tau_1+1}\,\sqrt{\varsigma+\Gamma_{k_{\tau_1}+1}}
\le \sqrt{\tau_1+1}\,\kappa_T.
\eeqn

Using the switching condition \req{switch} and the first part 
of \req{achle1}, we then deduce that, whether $k_\tau$ belongs
to $\{k_\nu\}$ or not,
\[
\sum_{\tau=\tau_0}^{\tau_1}\oN{k_\tau}
\le \sum_{\tau=\tau_0}^{\tau_1} \beta \aT{k_\tau}\oT{k_\tau}
\le \frac{\beta\eta}{\sqrt{\varsigma}} \sum_{\tau=\tau_0}^{\tau_1} \oT{k_\tau}
< \frac{\beta\eta\kappa_T\sqrt{\tau_1+1}}{\sqrt{\varsigma}}.
\]
Summing this bound with \req{tantan} then gives the second inequality
of \req{Tcomp}.  The first results from AS.8.
} 

\subsection{Normal complexity}

We now exploit the bound on $\oN{k_\nu}$ stated in Lemma~\ref{telescoping} to
analyze the complexity of the subsequence of normal iterations. We first show that the sum of the norms of constraint violations is bounded.

\llem{normal-complexity1}{
Suppose that AS.0--AS.9 hold. Then, for any $\nu_1>0$,
\beqn{Ncomp1}
\xi \sum_{\nu=\nu_0}^{\nu_1}\|c_{k_\nu}\|
\le \sum_{\nu=\nu_0}^{\nu_1}\oN{k_\nu}
< \kappa_N,
\eeqn
where
\beqn{kappaN-def}
\kappa_N = \kap{gap}+\kap{tan}\log\left(1+\frac{\kappa_T^2}{\varsigma}\right).
\eeqn
}

\proof{
The bound \req{full-telescopic} ensures that
\beqn{nc1-sto}
\sum_{\nu=\nu_0}^{\nu_1}\oN{k_\nu}
\leq \kap{gap}+\kap{tan}\log\left(1+\frac{\Gamma_{k_{\tau_1}+1}}{\varsigma}\right),
\eeqn
where $k_{\tau_1}$ is the index of the last tangential iteration before
$k_{\nu_1}$. Substituting the bound \req{Gamma-bound} in this inequality then gives
the second inequality of \req{Ncomp1}, the first resulting
again from AS.8.
}
\noindent
This allows us to derive boundedness of a combined primal and dual criticality measure.

\llem{normal-complexity}{
Suppose that AS.0--AS.9 hold. Then, for any $\nu_1\ge 0$,
\beqn{Ncomp}
\xi \sum_{\nu=\nu_0, k_\nu \not \in \{k_\tau\}}^{\nu_1}\Big(\oT{k_\nu}+\|c_{k_\nu}\|\Big)
\le\sum_{\nu=\nu_0, k_\nu \not \in \{k_\tau\}}^{\nu_1}\Big(\oT{k_\nu}+\oN{k_\nu}\Big)
< \kappa_N \left(1+\frac{\kappa_T}{\beta\eta}\right).
\eeqn
}

\proof{
Using the switching condition
\req{switch} for $k_\nu \not \in \{k_\tau\}$, we obtain that, for such
$k_\nu$ with $\nu \in \iibe{\nu_0}{\nu_1}$,
\beqn{fromswitch}
\oN{k_\nu} > \beta\,\aT{k_\nu}\oT{k_\nu}.
\eeqn
As in the previous lemma, let $k_{\tau_1}$ be the index of the last
tangential iteration before $k_{\nu_1}$.
Thus using \req{Gamma-bound}, 
\[
\alpha_{T,k_\nu}
= \frac{\eta}{\sqrt{\varsigma+\Gamma_{k_\nu}}}
\ge \frac{\eta}{\kappa_T}.
\]
Substituting this bound in \req{fromswitch}, we find that, for $\nu \in \iibe{\nu_0}{\nu_1}$,
\beqn{nc1b}
\oN{k_\nu}
\ge \frac{\beta\eta}{\kappa_T}\,\oT{k_\nu}.
\eeqn
With \req{Ncomp1}, this implies that
\[
\sum_{\nu=\nu_0}^{\nu_1}\oT{k_\nu}
\le \frac{\kappa_T}{\beta\eta}\sum_{\nu=\nu_0}^{\nu_1}\oN{k_\nu}
\le \frac{\kappa_N\,\kappa_T}{\beta\eta}.
\]
Summing this bound with \req{Ncomp1} and using AS.8 gives
\req{Ncomp}.
} 

\subsection{Combined complexity}

We finally assemble the pieces of the puzzle to derive our main result on the global rate of convergence of the \algname\ algorithm.

\lthm{final-complexity}{
Suppose that AS.0-AS.9 hold. Then, for any $k\ge 0$,
\beqn{the-complexity}
\frac{1}{k+1}\sum_{j=0}^k\Big(\oT{j}+\|c_j\|\Big)
\leq \frac{\kap{\algname,1}}{\sqrt{k+1}}+\frac{\kap{\algname,2}}{k+1}
= \calO\left(\frac{1}{\sqrt{k+1}}\right),
\eeqn
where
\[
\kap{\algname, 1}= \frac{\kappa_T}{\xi}\left(1+\frac{\beta\eta}{\sqrt{\varsigma}}\right)
\tim{ and }
\kap{\algname, 2}=\frac{\kappa_N}{\xi}\left(1+\frac{\kappa_T}{\beta\eta}\right).
\]
}

\proof{
Consider iterations of both types (tangential
and normal) from $0$ to $k$  by defining $\min[k_{\nu_0},k_{\tau_0}] = 0$ and
$\max[k_{\nu_1},k_{\tau_1}]= k$ (as in Lemma~\ref{telescoping}). We then
obtain, by combining \req{Tcomp} and \req{Ncomp}, that
\[
\begin{aligned}
\sum_{j=0}^k&\big(\oT{j}+\|c_j\|\Big)
= \sum_{\tau=\tau_0}^{\tau_1}\Big(\oT{k_\tau}+\|c_{k_\tau}\|\Big)
+\sum_{\nu=\nu_0, k_\nu \not \in \{k_\tau\}}^{\nu_1}\Big(\oT{k_\nu}+\|c_{k_\nu}\|\Big)\\
&\le \frac{\kappa_T}{\xi}\sqrt{k+1}\left(1+\frac{\beta\eta}{\sqrt{\varsigma}}\right)
+ \frac{\kappa_N}{\xi}\left(1+\frac{\kappa_T}{\beta\eta}\right),
\end{aligned}
\]
where we used the inequalities  $\tau_1 \le k_{\tau_1}\leq k$ and
$k_{\nu_1}\leq k$. The bound \req{the-complexity} is finally obtained
by dividing both sides by $k+1$.
}

\noindent
Remarkably, Theorem~\ref{final-complexity} implies that obtaining an $\epsilon$-approximate first-order critical point, that is an iterate $x_k$ such that $\oT{j}+\|c_j\| \le \epsilon$, requires at most $\calO(\epsilon^{-2})$ iterations of the \al{ADIC} algorithm, a complexity which is, in order, the same as that of steepest-descent and Newton's methods on unconstrained problems \cite[Theorems~2.2.2 and 3.1.1]{CartGoulToin22}.

\section{Numerical illustration}\label{sec4}

We now illustrate the behaviour of three variants of the \algname\ algorithm on problems from the {\sf CUTEst} \cite{GoulOrbaToin15b} collection as provided in Matlab by S2MPJ \cite{GraToin25}. All nonlinear optimization problems in the collection involving general constraints and at most 200 variables were considered, leading to a test set of 312 problems. The algorithmic variants are defined as follows.
\begin{itemize}
\item The first variant (ADIC-LP) follows Section~\ref{ss:LP} and computes the dual criticality measure and tangential step using the linear optimization subproblems \req{LP-chiT} and \req{LP-sT}, respectively.
\item The second variant (ADIC-BK) again follows Section~\ref{ss:LP} and computes the dual criticality measure using \req{LP-chiT}, but then uses the simple formula \req{sT-backtrack} to define the tangential step.
\item The third variant (ADIC-PR) uses the projection approach of Section~\ref{ss:proj}, in which the dual criticality is given by \req{pik-def} and the tangential step is defined by \req{proj1-sT}.
\end{itemize}
All three variants have been (trivially) extended to handle general lower and upper bounds on the variables (instead of mere non-negativity constraints), thereby making them applicable to general constrained problems (after transformation of inequality constraints into equalities and the introduction of slack variables, if needed).

Because the variants use different criticality measures, a uniform (external) termination criterion was implemented
in order to enforce consistency in the comparison. For all variants, a problem was considered solved as soon as 
\[
\chi_{T,k} \le 10^{-4} \tim{ and }\chi_{N,k} \le 10^{-5},
\]
where $\chi_{T,k}$ and $\chi_{N,k}$ are defined in \req{chiTk-def} and \req{chiNk-def}, respectively. Note that this accomodates the (unfortunate but unavoidable) case where an infeasible minimizer of the equality constraint's violation is found (the bound constraints are satisfied throughout the algorithms). This situation is excluded from our theoretical analysis by AS.8 but does occur in practice. A maximum number of 50000 iterations and a 3600 seconds time limit were imposed. Finally, the  algorithmic parameters were chosen as
\[
\varsigma = 10^{-5}, \ms \eta = 2, \ms \theta_T = 1, \ms \theta_N = 5, \ms\kappa_n = 10^{-2} \tim{and} \beta = 10^3.
\]

\begin{figure}[htb]
\includegraphics[width= 0.48\linewidth]{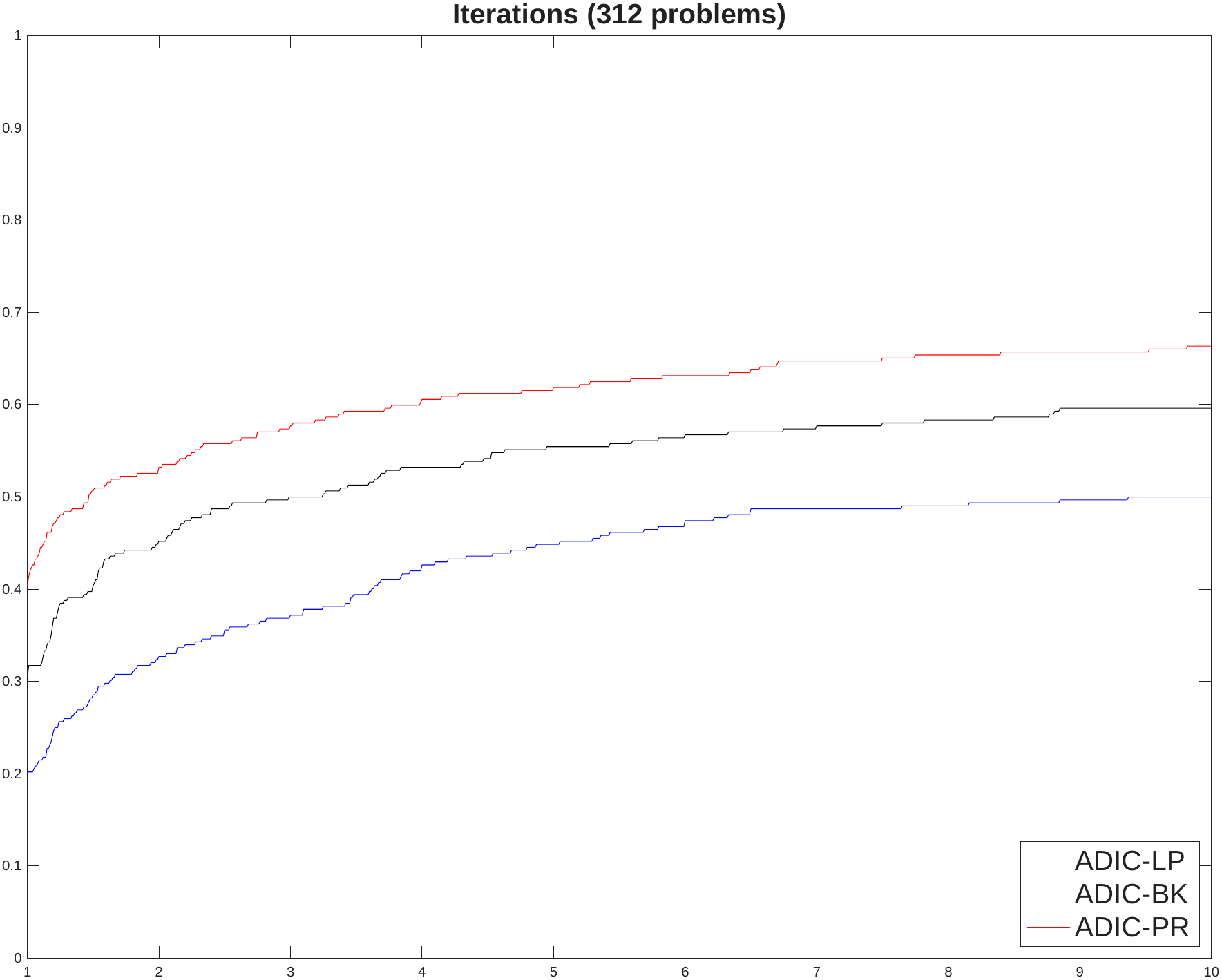}
\includegraphics[width= 0.48\linewidth]{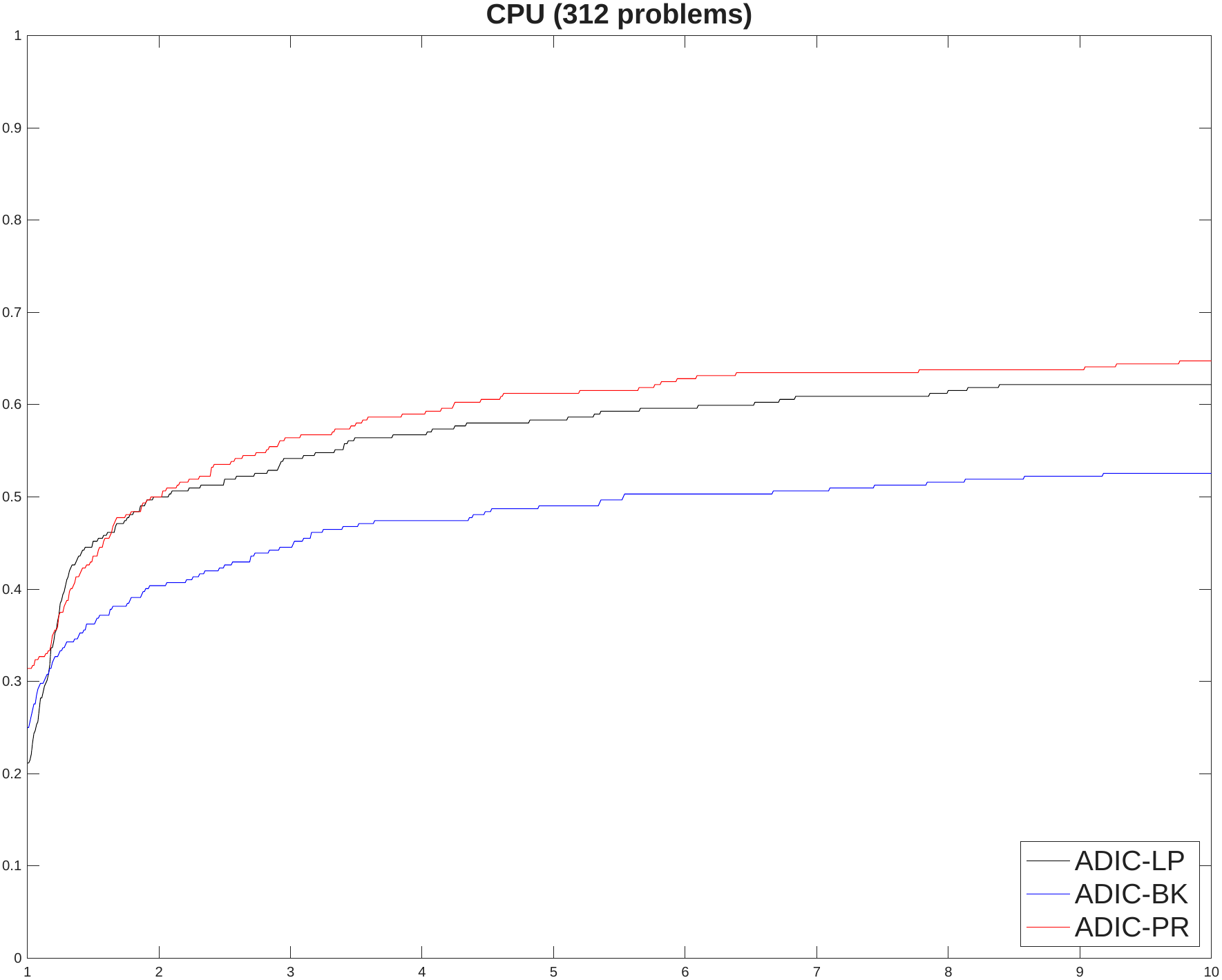}
\caption{\label{fig:noiseless}Iteration (left) and CPU time (right) performance profiles for three variants of \algname\ on 312 constrained {\sf CUTEst} problems (noiseless gradients)}
\end{figure}

We first report on a set of experiments in which the gradients used by the three variants were exact. To summarize the results, we computed three performance statistics: efficiency in terms of iterations, efficiency in terms of CPU time needed and reliability. The latter, which we denote by "Rel" in what follows, is simply computed as the percentage of successfully solved problems. For the two first efficiency statistics, we follow the approach of \cite{PorcToin17} and compute, for each variant, the area below the relevant curve in a performance profile comparing the three variants, truncated at a ``ratio to best performance'' equal to 10. The iteration-based statistic is denoted by ``Iters'' and the CPU-based one by ``Time''. Values of these statistics should be as close to one as possible. Results are presented in Table~\ref{tab:noiseless}.  The corresponding iteration and CPU performance profiles are shown in Figure~\ref{fig:noiseless}.

\begin{table}[htb]
\begin{center}
\begin{tabular}{|l|c|c|c|}
\hline
Variant  &    Iters  &      Time    &     Rel   \\
\hline
ADIC-LP  &  0.54  &   0.57   &   68.27 \\
ADIC-BK  &  0.43  &   0.48   &   61.54 \\
ADIC-PR  &  0.61  &   0.59   &   71.15 \\
\hline
\end{tabular}
\caption{\label{tab:noiseless}Efficiency and reliability statistics for three variants of \algname\ on 312 constrained {\sf CUTEst} problems (noiseless gradients)}
\end{center}
\end{table}

Table~\ref{tab:noiseless} and Figure~\ref{fig:noiseless} indicate that the projection-based ADIC-PR clearly outperforms both ADIC-LP and ADIC-BK, on all 3 statistics. The dominance of ADIC-PR over ADIC-LP in CPU time is however marginal, despite the fact that two linear optimization subproblems must be solved at each iteration of ADIC-LP, against a single projection subproblem for ADIC-PR. The ADIC-BK variant, which only requires the solution of a single linear optimization problem per iteration, remains slower mostly because it typically needs more iterations than other variants per successfully solved problem.

We finally show that our claim that OFFO methods are reliable in the presence of noise is vindicated in practice.  To analyze this, we considered the subset of our problems by keeping those that were solved in the absence of noise by at least one variant, giving a set of 247 test problems. We then added relative random Gaussian\footnote{With zero mean and unit variance.} noise of increasing magnitude (5\%, 15\%, 25\% and 50\%) to the gradients of the objective function and ran each problem 20 times independently, with $\chi_{T,k} \le 10^{-3}$ and $\chi_{N,k} \le 10^{-3}$. We then computed the total reliability of our three variants on the resulting 4940 runs for each of the 5 noise levels.
The results are presented in Table~\ref{tab:noisy}. They show an impressive stability for increasing noise levels, and indicate that, in our view remarkably, \algname\ is capable of handling very substantial perturbations of the gradient of the objective function (50\% relative noise results in barely one significant digit in the gradient) for a reasonable accuracy requirement. It is also interesting to note that 
ADIC-LP becomes marginally more reliable than ADIC-PR for large noise levels.

\begin{table}
    \begin{center}
        \begin{tabular}{|l|ccccc|}
        \hline
        Variant & 0\% & 5\% & 15\% & 25\% & 50\% \\
        \hline
        ADIC-LP & 81.78 & 74.55 & 72.85 & 73.14 & 72.73 \\
        ADIC-BK & 74.09 & 69.11 & 67.11 & 65.77 & 59.37 \\
        ADIC-PR & 89.88 & 77.79 & 74.82 & 72.57 & 67.19 \\
        \hline
        \end{tabular}
         \caption{\label{tab:noisy} Reliability statistics for three variants of \algname\ on 247 constrained {\sf CUTEst} problems for relative random Gaussian noise levels of 0\%, 5\%, 15\%, 25\% and 50\% on the objective function's gradient}
    \end{center}
\end{table}

\section{Conclusions and perspectives}\label{sec5}

We have proposed a new OFFO algorithm for the solution of smooth optimization problems, with excellent stability in the presence of noise on the objective function's gradient. This "trust-funnel" algorithm uses adaptive switching between a normal step (reducing constraint violation), and tangential steps (improving dual optimality), the latter being inspired by the AdaGrad-norm algorithm \cite{DuchHazaSing11,ward2020adagrad} for unconstrained problems.  We have also provided a full analysis of the method's worst-case iteration complexity, showing that its global rate of convergence is, for problems with full-rank Jacobians, identical in order to the (optimal) rate of steepest-descent and Newton's method on unconstrained problems. This also provides an evaluation complexity for evaluations of the objective function's gradient, because each iteration requires a single gradient computation. Evaluation complexity for the constraint function and Jacobian is not direct and depends on the algorithm used in the normal step. We have finally conducted illustrative numerical experiments suggesting that the algorithm's performance and reliability are satisfactory (although admittedly not state-of-the-art) on noiseless problems, but that its reliability in the presence of significant noise on the objective function's gradient is very remarkable. 

Many questions remain for further investigation, including the incorporation of second-order information, should it be available, a component-wise version of the algorithm (closer to AdaGrad as opposed to AdaGrad-norm) and a full stochastic complexity analysis. These topics are the subject of ongoing research. Exploiting the independent structure of normal and tangential steps to allow for specific preconditioning of the normal step and relaxing the full-rank assumption on the Jacobians are also of interest.

{\footnotesize
  
\section*{\footnotesize Acknowledgement}

Philippe Toint is grateful for the continued and
friendly support of the APO team at Toulouse IRIT (F) and of DIEF at the
University of Florence (I).

\appendix
\appnumsection{Proof of Lemma~\ref{Lipschitz-things}}

We prove the five statements Lemma~\ref{Lipschitz-things} for
arbitrary $x,y \ge 0$.
\begin{enumerate}
\item That $c(x)$ is Lipschitz continuous with constant $L_c = \kappa_J$ directly follows from AS.4.
\item We have, from AS.3, AS.4 and AS.7 that
\[
\begin{aligned}
\|J(x)^Tc(x)-J(y)^Tc(y)\|
&= \left\|\Big(J(x)-J(y)\Big)^Tc(x)+J(y)^T\Big(c(x)-c(y)\Big)\right\|\\
&\le (\kappa_c L_J + \kappa_J L_c) \|x-y\|\\
&= (\kappa_c L_J + \kappa_J^2)\|x-y\|,
\end{aligned}
\]
yielding $L_{JTc} = \max[1,\kappa_c L_J + \kappa_J^2]$.
\item Define $A(x) = J(x)J(x)^T$. AS.5 then implies that $A(x)$ is (symmetric) positive-definite with smallest eigenvalue bounded below by $\sigma_0^2$. As a consequence, $\widehat\lambda(x)$ is well defined by \req{eq:ls-mult}.
\item Moreover, AS.2 and AS.4 then imply that
\[
\|\widehat\lambda(x)\| \le \frac{\kappa_g\kappa_J}{\sigma_0^2},
\]
yielding $\kappa_\lambda=\kappa_g\kappa_J/\sigma_0^2$. We have also,
using AS.4 and AS.7, that
\[
\|A(x)-A(y)\|
\leq \left\|\big(J(x)-J(y)\Big)J(x)^T + J(y)\Big(J(x)-J(y)\Big)^T\right\|
\le 2 \kappa_JL_J \|x-y\|.
\]
from which we deduce that
\[
\|A(x)^{-1} -A(y)^{-1}\|
= \left\|A(x)^{-1}\Big(A(x)-A(y)\Big)A(y)^{-1}\right\|
\le \frac{2\kappa_J L_J}{\sigma_0^4}\|x-y\|.
\]
We also have that
\[
\|J(x)g(x)-J(y)g(y)\|
= \left\|\Big(J(x)-J(y)\Big)g(x) + J(y)\Big(g(x)-g(y)\Big)\right\|
\le (\kappa_g L_J+\kappa_J L_g)\|x-y\|.
\]
where we used AS.2, AS.4, AS.6 and AS.7.
Therefore, using \req{eq:ls-mult},
\[
\begin{aligned}
\|\widehat\lambda(x)-\widehat\lambda(y)\|
&= \left\|\Big(A(x)^{-1} -A(y)^{-1}\Big)J(x)g(x)+A(y)^{-1}\Big(J(x)g(x)-J(y)g(y)\Big)\right\|\\
&\le \frac{1}{\sigma_0^2}\left(\frac{2\kappa_g \kappa_J^2 L_J}{\sigma_0^2}+ \kappa_g L_J+\kappa_JL_g\right)\|x-y\|,
\end{aligned}
\]
yielding $L_\lambda = \Big((2\kappa_g \kappa_J^2 L_J)/\sigma_0^2+ \kappa_g L_J+\kappa_JL_g\Big)/\sigma_0^2$.
\item Finally, we obtain that, for any $\lambda$ such that $\|\lambda\|\le\kappa_\lambda$,
\[
\|\nabla_xL(x,\lambda)-\nabla_yL(y,\lambda)\|
= \left\|g(x)-g(y)+\Big(J(x)-J(y)\big)^T\lambda \right\|
\le \left(L_g+\kappa_\lambda L_J\right)\|x-y\|,
\]
yielding $L_L = L_g+\kappa_\lambda L_J= L_g+\kappa_g\kappa_J L_J/\sigma_0^2$.
\end{enumerate}
\end{document}